\documentclass[12pt]{article}

\usepackage{tikz}
\usepackage{graphicx}  
\usepackage{float}
\usepackage{amsmath,amssymb,amsfonts,theorem,makeidx,latexsym,epsfig,subfigure}

\usepackage{color}

\newtheorem{defn}{Definition}[section]

\newtheorem{lemma}[defn]{Lemma}

{\theorembodyfont{\rmfamily}

	\newtheorem{ex}[defn]{Example}}

\newtheorem{thm}[defn]{Theorem}

\newtheorem{prop}[defn]{Proposition}

\newtheorem{cor}[defn]{Corollary}

\newtheorem{rem}[defn]{Remark}

\numberwithin{equation}{section}

\newcommand{\h}{{\cal H}}

\newcommand{\ltn}{{\ell}^2(\mathbb N)}

\newcommand{\si}{S^{-1}}

\newcommand{\mn}{\mathbb N}

\newcommand{\mr}{\mathbb R}

\newcommand{\mc}{\mathbb C}

\newcommand{\mcn}{{\mathbb C}^n}

\def\bp{{\noindent\bf Proof. \ }}

\def\ep{\hfill$\square$\par\bigskip}

\def\bqs{\begin{equation}}

	\def\eqs{\tag*{$\square$}\end{equation}\par\bigskip}

\def\la{\langle}

\def\ra{\rangle}

\def\ftn{\{f_k\}_{k=1}^n}

\def\ctk{\{c_k\}_{k=1}^\infty}

\def\etk{\{e_k\}_{k=1}^\infty}

\def\ftkn{\{f_k\}_{k=1}^n}

\def\suk{\sum_{k=1}^\infty}

\def\vn{\vspace{.1in}\noindent}

\def\bop{\begin{op}\rm}

	\def\eop{\end{op}}

\def\bee{\begin{eqnarray}}

	\def\ene{\end{eqnarray}}

\def\bes{\begin{eqnarray*}}

	\def\ens{\end{eqnarray*}}

\def\bei{\begin{itemize}}

	\def\eni{\end{itemize}}

\def\bt{\begin{thm}}

	\def\et{\end{thm}}

\def\bc{\begin{cor}}

	\def\ec{\end{cor}}

\def\bpr{\begin{prop}}

	\def\epr{\end{prop}}

\def\bl{\begin{lemma}}

	\def\el{\end{lemma}}

\def\bd{\begin{defn}}

	\def\ed{\end{defn}}

\def\bex{\begin{ex}}

	\def\enx{\end{ex}}

\def\bfi{\begin{fig}}

	\def\efi{\end{fig}}

\newcommand{\C}{\mathbb C}

\newcommand{\cd}{\C^{d}}

\newcommand{\fn}{\{f_k\}_{k=1}^n}
\newcommand{\fd}{\{f_k\}_{k=1}^d}

\newcommand{\tnf}{\{T^{k-1}f_1\}_{k=1}^n}
\newcommand{\tdf}{\{T^{k-1}f_1\}_{k=1}^d}

\def\gtkd{\{g_k\}_{k=1}^d}
\def\ftkd{\{f_k\}_{k=1}^d}

\def\tng{\{T^{k-1}g_1\}_{k=1}^n}
\def\tdg{\{T^{k-1}g_1\}_{k=1}^d}

\def\br{\begin{rem}}

	\def\er{\end{rem}}

\title{Cyclic frames in finite-dimensional Hilbert spaces}

\date{\today}

\author{}

\begin{document}
	
	\author{Ole Christensen, Navneet Redhu, Niraj K. Shukla}
	
	\maketitle
	
	\begin{abstract}  Generalizing a definition by Kalra \cite{Kalra}, the purpose of this paper is to analyze cyclic frames in finite-dimensional Hilbert spaces. Cyclic
		frames form a subclass of the dynamical frames introduced and analyzed in detail by Aldroubi et al. in \cite{ACM} and subsequent papers; they are particularly interesting due to their
		attractive properties in the context of erasure problems. By applying an alternative approach, we are able to shed new light on general dynamical frames as well as cyclic frames. In particular, we provide a characterization of dynamical frames, which in turn leads to
		a characterization of cyclic frames. 
	\end{abstract}
	
	\begin{minipage}{120mm}
		
		
	\end{minipage}
	\

	\section{Introduction}
	
	Dynamical sampling was introduced around 2013--2014 by Aldroubi et al. in a series of papers \cite{AJI,ACM, ACCMP}. Discarding the applied context in which dynamical sampling was introduced, the mathematical issue is how to construct frames in a Hilbert space by iterating the action of a bounded operator on a collection of vectors; in the literature such frames are called  {\it dynamical frames.} Very fast, dynamical sampling became an active research area, mainly with contributions dealing with infinite-dimensional spaces; note, however, that
	already the paper \cite{ACM} characterized the frame property for such iterated systems in the finite-dimensional setting.  Other theoretical contributions are contained in the paper \cite{AJPAM} by Ashbrock and Powell.
	
	In the current paper we focus on dynamical frames in finite-dimensional
	spaces; for convenience, and without loss of generality, we formulate all results for 
	frames in $\cd.$  We will complement the analysis in \cite{ACM} with an alternative approach, which shed light on dynamical frames from a different angle. Based on this approach, we will analyze
	{\it cyclic frames,} a concept that is strictly restricted to the finite-dimensional setting.  Cyclic frames were introduced in the paper \cite{Kalra} by Kalra in 2006; since the terminology has not been used in later papers, we will allow ourself to use the same name for a slightly more general concept, see Section   \ref{251302b} for details. Already Kalra noticed that cyclic frames have attractive features in the context of erasure problems. For more details on erasures, we refer to \cite{HRBPVI,Kalra,LD,LD2}.
	The main purpose of our paper is to characterize (our generalized version of) cyclic frames and  discuss their key features.
	
	Our alternative approach to dynamical frames is presented in Section  \ref{251302a}. In Section \ref{251302b} we state the formal definition of cyclic frames
	and two ways  of characterization of such frames.  Furthermore we prove a number
	of properties that follow directly from our approach in Section \ref{251302a}.  In Section \ref{251302c} we consider tight cyclic frames and their connection to erasure problems. Finally, the appendix contains proofs of a number of more technical results.

	\section{Dynamical frames} \label{251302a}

	A frame $\fn$ for $\cd$ is called a \textit{dynamical frame} for $\cd$  if there exists a linear operator $T: \cd \to \cd$   such that
	\begin{equation}\label{eq1}
		f_{k}=T^{k-1} f_1 \mbox{ for } 1 \leq k \leq n.
	\end{equation} 
	In that case we can write  \bee \label{252103a} \fn=\tnf= \{f_1, Tf_1, \dots, T^{n-1}f_1\}.\ene 
	Note that in \cite{AJPAM}, a dynamical frames is defined as a frame of the form
	\bee \label{252103b}
	\fn= \{T^{k}f\}_{k=1}^n =
	\{Tf, T^2f, \dots, T^{n-1}f\}, f \in \cd;\ene clearly, this is a subclass of the dynamical frames considered in the current paper.

	 Our main purpose is to
	consider a special class of dynamical frames, to be introduced in Section \ref{251302b}. For this purpose we first state a number of properties of general dynamical frames. The following lemma is well-known; it implies in particular that every basis for $\cd$ is a dynamical frame.
	
	\bl \label{lem_linearly_independent_dynamical_representation} 
	Assume that $\fn$ is a collection of linearly independent vectors in  $\cd.$  Then there exists a linear operator $T : \cd \to \cd$ such that 
	\bee \label{eq_fk=Tnf} 
	\fn	 = \tnf.
	\ene  
	\el
	\bp
	The assumption of linear independence of $\fn$ implies that $n \leq d$. 
	Due to linear independence of the vectors $\fn,$ we can clearly  define the operator $T$ by $Tf_k = f_{k+1}$ for $k = 1, \dots, n - 1.$   Note that this only defines $T$ on $\text{span}\{f_1, \dots, f_{n-1}\}$. The operator can be extended to a linear operator on $\cd$ by defining it in an arbitrary way on a basis for $\text{span}\{f_1, \dots, f_{n-1}\}^\perp$.
	\ep
	
	In general, the key feature of a frame is that it might contain more elements than a basis. 
	The following result states that if $\tnf$ is a dynamical frame for $\cd$, then the set of the {\it first} $d$ consecutive vectors, i.e., the set $\tdf$, forms a basis for $\cd$. The similar result for dynamical frames of the form \eqref{252103b} was proved in \cite{AJPAM}.
	
	\bl \label{first_d_vector_basis}
	Assume that $\tnf$ is a frame for $\cd$. Then the first $d$ elements, i.e., $\tdf$, form a basis for $\cd$.
	\el
	
	\bp
	Assume that $\tdf$ is not a basis for $\cd$. Then $\tdf$ is linearly dependent. There are now two possibilities:
	
	1) If $\{T^{k-1}f_1\}_{k=1}^{d-1}$ are linearly independent, then $T^{d-1}f_1 \in \text{span}(\{T^{k-1}f_1\}_{k=1}^{d-1})$. But then the subspace $\text{span}(\{T^{k-1}f_1\}_{k=1}^{d-1})$ is invariant under the operator $T$, and, regardless of the choice of $n \in \mathbb{N}$, a family $\tnf$ cannot be a frame for $\cd$. This contradicts the assumptions in the lemma.
	
	2) If $\{T^{k-1}f_1\}_{k=1}^{d-1}$ are linearly dependent, the same argument as in 1) implies that $\text{span}(\{T^{k-1}f_1\}_{k=1}^{d-1})$ is invariant under $T$. Hence, again the subspace $\text{span}(\{T^{k-1}f_1\}_{k=1}^{d-1})$ is invariant under the operator $T$, and, regardless of the choice of $n \in \mathbb{N}$, a family $\tnf$ cannot be a frame for $\cd$. This contradicts the assumptions in the lemma.
	
	Altogether, we conclude that $\tdf$ is a basis for $\cd$.
	\ep
	
	Note that Lemma \ref{first_d_vector_basis} does not guarantee that any other set of $d$-consecutive vectors in $\tnf$ will form a basis for $\cd$; see Corollary \ref{251302e} for a more detailed discussion. 
	
	Our approach in the current paper is based
	on the next result, which  shows that the class of all finite dynamical frames is parametrized by a choice of $d$ linearly independent vectors $\fd$, a vector $\varphi \in \cd$ and an integer $n \ge d$.

	\bt \, \label{thm_independent_vector_dynamical_frame}
	
	\begin{itemize}
		\item[(i)] Consider a collection of linearly independent vectors $\fd$ in $\cd$, an arbitrary vector $\varphi \in \cd$, and an integer $n \geq d$. Define the linear map $T: \cd \to \cd$ by
		\begin{equation}\label{eq_Tf_kp}
			T f_k = f_{k+1}, \; k = 1, \ldots, d-1, \; T f_d = \varphi.
		\end{equation}
		Then $\tnf$ is a frame for $\cd$.
		\item[(ii)] Conversely, any  dynamical frame $\tnf$ for $\cd$ corresponds to the setup in (i), with $f_k = T^{k-1} f_1$, $k = 2, \ldots, d$, $\varphi = T f_d$.
	\end{itemize}
	\et
	
	\bp The statement in (i)  is precisely the construction in Lemma \ref{lem_linearly_independent_dynamical_representation}; already in the proof of the lemma it was highlighted that the construction works for arbitrary choices of the vector $\varphi \in \cd$. 
	Clearly $\fd= \{T^{k-1}f_1\}_{k=1}^d$ is a basis for $\cd,$ and hence $\tnf$ is an 
	(overcomplete) frame. 	
	Concerning (ii), Lemma \ref{first_d_vector_basis} shows that if $\tnf$ is a dynamical frame, then the first $d$ elements, $\tdf$, form a basis for $\cd$, i.e., a linearly independent set of vectors. Now the result follows from Lemma \ref{lem_linearly_independent_dynamical_representation}.
	\ep
	
	The next result concerns the special case of
	$d+1$ vectors in $\cd.$ We return to this setting in Lemma \ref{251302g} and Proposition \ref{252003a}.

	\bc A set of $d+1$ vectors $\{ f_k\}_{k=1}^{d+1}$ is a dynamical frame for $\cd$ if and only if $\{ f_k\}_{k=1}^{d}$ is a basis for $\cd$.
	\ec
	\bp Suppose $\{f_1,\ldots, f_{d+1}\}$ is a dynamical frame for $\cd$. Then by Lemma   \ref{first_d_vector_basis}, the first $d$ elements, i.e., \(\fd\), forms a basis for \(\cd\). The converse follows from Theorem \ref{thm_independent_vector_dynamical_frame}(i) by choosing \(n=d+1\) and \(\varphi=f_{d+1}\).
	\ep

	Note that in the setup of Theorem \ref{thm_independent_vector_dynamical_frame}(i), the operator $T$ in \eqref{eq_Tf_kp} is either surjective (if $\varphi \notin \text{span}\{f_2, \dots, f_d\}$), or its range has codimension one (if $\varphi \in \text{span}\{f_2, \dots, f_d\}$). The resulting frame has different properties in these two cases:
	
	\bc \label{251302e}
	\textit{Consider the setup in Theorem \ref{thm_independent_vector_dynamical_frame}(i). Then the following hold:}
	\begin{itemize}
		\item[(i)] If the operator $T$ in \eqref{eq_Tf_kp} is surjective, each collection of $d$ consecutive vectors $\{T^{k-1}f_1\}_{k=\ell+1}^{\ell+d}$ is a basis for $\cd$, for any $\ell =0, 1, \dots$.
		\item[(ii)] If the operator $T$ in \eqref{eq_Tf_kp} is not surjective, $\{T^{k-1}f_1\}_{k=\ell+1}^{\ell+d}$ is not a basis for $\cd$, for any $\ell = 1, 2\dots$.
	\end{itemize}
	\ec
	
	\bp For the proof of (i), given any $\ell \in \{0,1,2,\dots\}$, the vectors $\{T^{k-1}f_1\}_{k=\ell+1}^{\ell+d}$ are the images of the basis vectors $\tdf$ under the bijective operator $T^{\ell}$, and hence itself a basis for $\cd$.  
	Concerning (ii), if $T$ is not surjective, the vectors $\{T^{k-1}f_1\}_{k=\ell+1}^{\ell+d}$ cannot span $\cd$, and hence not be a basis (or frame) for $\cd$.
	\ep
	
	Note that for dynamical frames of the form \eqref{252103b}, it was already observed in \cite{AJPAM} that every collection of $d$ consequtive vectors form a basis for $\cd.$ 
	
	Let us end this section with an application of Theorem \ref{thm_independent_vector_dynamical_frame}. Using a measure-theoretic approach, Ashbrock and Powell showed in \cite{AJPAM} that every frame
	$\fn$ (dynamical or not)  for $\cd$ has a dynamical dual frame; see \cite[Theorem 3.6]{AJDMRJBMIPAM} for an alternative approach. If we
	allow ourself to reorder the elements in the frame $\fn,$ we can provide a short and constructive proof of this:
	
	\bpr \label{251103a} Let $\fn$ be a frame for $\cd,$ ordered such that $\ftkd$ is linearly independent. Then $\fn$ has a dual frame of the form
	$\tng.$ \epr
	
	\bp  Every frame $\fn$ contains a basis, so by a reordering we can
	assume that the first $d$ elements  $\ftkd$ is a basis for
	$\cd.$ Let $\gtkd$ denote the dual basis, and define according to
	Theorem \ref{thm_independent_vector_dynamical_frame}
	the operator $T: \cd \to \cd$ by
	\bes Tg_k:= g_{k+1}, \, k=1, \dots, d-1, \, Tg_k:=0.\ens  Then $\gtkd= \tdg;$ furthermore,
	by construction, 
	\bes \tng= \{g_1, g_2, \dots, g_d, 0,0, \dots 0\},\ens which is clearly a dual frame
	of $\ftkn.$ \ep

	\section{Cyclic frames} \label{251302b}
	We now move to the cyclic frames. The term {\it cyclic frame}
	was coined in the paper \cite{Kalra} by Kalra in 2006.  The
	terminology has apparently not been used in the literature since 
	this first paper;  we will use the name in a more general sense, as follows.

	\bd \label{251302f}  A dynamical frame $\fn=\tnf$ for $\cd$ is called a cyclic frame if $T^n=I;$ if $n$ is the minimal choice of a positive integer such that $T^n=I,$ we call
	$\tnf$ a minimal cyclic frame.  \ed 
	
	The difference between the terminology in \cite{Kalra} and
	Definition \ref{251302f} is that Kalra requires the operator $T$ to be unitary; furthermore, the distinction between cyclic frames and minimal cyclic frames does not occur in \cite{Kalra}.  
	
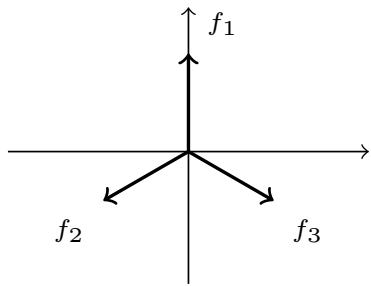
\begin{figure}[H] \label{252103d}
		\centering
		\scalebox{1.6}{ 
			\begin{tikzpicture}
				\draw[->] (-1.5,0) -- (1.5,0) node[below] {};
				\draw[->] (0,-1.1) -- (0,1.2) node[left] {};
				
				\draw[->,thick,black] (0,0) -- (0,0.8164965809) node[above right] {\tiny $f_1$};
				\draw[->,thick,black] (0,0) -- (-0.7071067812,-0.4082482905) node[below left] {\tiny $f_2$};
				\draw[->,thick,black] (0,0) -- (0.7071067812,-0.4082482905) node[below right] {\tiny $f_3$};

				\node at (0,0) [below left] {};
			\end{tikzpicture}
		}
		\caption{Mercedes-Benz frame $\{f_k\}_{k=1}^3$ in $\mr^2$.}
		\label{fig:Mercedes-Benz frame}
	\end{figure}

	Cyclicity of a frame $\fn=\tnf$ means that if we consider iterates
	$\{ T^{k-1}f_1\}_{k=1}^N$ for some $N>n,$ the resulting system of vectors would simply repeat (some of) the vectors in $\fn;$ and $\tnf$ being
	minimal means that there are no repetitions at all.
	In order to have an intuitive feeling for this (and for the subsequent results), we recommend the reader to think about the so-called {\it Mercedes-Benz frame}
	$\{f_k\}_{k=1}^3$ for $\mr^2,$  which is a minimal cyclic frame with $T$ being a rotation of $\frac{2\pi}{3}$ rad; see Figure \ref{fig:Mercedes-Benz frame}. We return
	to this frame in Example \ref{252103c}.
	The approach in Section \ref{251302a} has a number of immediate consequences for cyclic frames:

	\bc A basis $\fd$ for $\cd$  is a cyclic  frame. 
	\ec
	
	\bp This follows from  Theorem  \ref{thm_independent_vector_dynamical_frame}(i), with 
	$n=d$ and $\varphi=f_1.$     \ep 
	
	A cyclic frame for $\cd$ has the particular property that each collection of $d$ consecutive elements form a
	basis for $\cd;$ as we saw in Corollary \ref{251302e}(ii) general dynamical frames do not have this property.
	
	\bc\label{cor_d_consucutive_vectors_forms_basis}
	Let  $\tnf$ be a cyclic  frame for  $\cd$, with $n>d$.  Then, for every $\ell \in \{1,\cdots, n-d\}$, the set $\{T^{n-1} f_1\}_{n= \ell+1}^{\ell +d}$ is a basis for $\cd$.
	\ec
	\bp Since $T$ is a cyclic operator,  $T$ is  invertible, and
	hence  surjective. Now the result follows from Corollary \ref{251302e}(i).
	\ep
	
	We will now address the question of how to construct cyclic frames $\tnf$ for $\cd,$ with $n>d.$  For this purpose we refer to Theorem \ref{thm_independent_vector_dynamical_frame}(i), which shows that
	the construction of {\it any} dynamical frame $\fn =\tnf, \, n>d,$
	must be based on a choice of a basis 
	$\fd \subset \cd$ and a vector $\varphi \in \cd;$ the linear map $T: \cd \to \cd$ must then be defined by
	\begin{equation}\label{eq_Tf_k}
		T f_k = f_{k+1}, \; k = 1, \ldots, d-1, \; T f_d = \varphi.
	\end{equation}

	Writing the vector $\varphi \in \cd$  as
	\[
	\varphi = c_1 f_1 + c_2 f_2 + \cdots + c_d f_d, \, \, c_1, \dots, c_d\in \mc,
	\]
	we see that to construct  a cyclic frame is a question of choosing $c_1, \dots, c_d$ such that $T^n = I$ for some $n >d$. This can easily be done in the case $n=d+1;$ we
	return to the construction below in Proposition \ref{252003a}.

	\bl \label{251302g}
	Let \(\fd\) be a basis for \(\cd\). Then, with \(f_{d+1} := -f_1 - f_2- \cdots - f_d\),
	the family \(\{f_k\}_{k=1}^{d+1}\) is a cyclic frame.
	\el
	
	\bp
	Following the construction in Theorem \ref{thm_independent_vector_dynamical_frame}(i), and taking \(n = d+1\) 
	and \(\varphi := f_{d+1} = -f_1 - f_2- \cdots - f_d\), we have that
	\bes
	T f_{d+1} &=& -T f_1 - T f_2 - \cdots - T f_d\\
	&=& - f_2 - f_3 - \cdots - f_d - (- f_1 - f_2 - \cdots - f_d) = f_1.
	\ens
 From the definition of $T$ in \eqref{eq_Tf_kp}, we have $f_{d+1} = T^{d} f_1$.  Therefore,
	\bee\label{eq_revis_3.2}
	T^{d+1} f_1 = T f_{d+1} = f_1.
	\ene
	Applying $T^{k-1}$ to both sides of \eqref{eq_revis_3.2} for each $k \in \{1,2,\dots,d\}$, and using the definition of $T$ from \eqref{eq_Tf_kp}, we obtain
	\bes
	T^{d+1} f_k = f_k, \quad \forall\, k \in \{1,2,\dots,d\}.
	\ens
	Since $\fd$ forms a basis of $\cd$, it follows that $T^{d+1} = I$.  Hence, $\{f_k\}_{k=1}^{d+1}$ is a cyclic frame.
	\ep
	
	Notice that for $n=2$ and $f_1,f_2$ being the two first vectors
	in the Mercedes-Benz frame, the construction in Lemma \ref{251302g} yields precisely the third vector $f_3.$ 
	It is clear that by looking at the $n$th roots of unity in $\mc,$
	we can easily construct cyclic frames for $\mr^2$ (and hereby for
	$\mc^2$) for an arbitrary $n>2.$ We will now state a construction of a cyclic frame $\tnf$ for $\cd,$ for {\it any} choice of $d\in \mn $ and $n>d.$ The basic insight in Lemma \ref{lem_construction_of_unitary_cyclic_dynamical_frames} (iii) below goes back to Zimmermann \cite{G}, who proved the result for the case
	where $f_1=(1, \dots, 1);$ also, Kalra \cite{Kalra} noticed that the obtained frame is cyclic. The result is based on the $n$th roots of unity, which we write in the form
	\bee \label{252802a} \omega= e^{\frac{2\pi i m}{n}}, \, m=1, \dots, n.\ene Recall that if $gcd(m,n)=1,$  $\omega$ is said to be a {\it primitive}
	root of order $n;$  this means that while $\omega^n=1,$ we also have that $\omega^k \neq 1$ for $k=1, \dots, n-1.$

	\bl \label{lem_construction_of_unitary_cyclic_dynamical_frames}
	Fix any \(n > d\), and let $\omega_1, \omega_2, \dots, \omega_d$ denote  $n$th roots of unity. Consider the \(d \times d\) matrix
	\[	T := 
	\begin{pmatrix}
		\omega_1 & 0 & 0 & \cdot &\cdot & 0 \\
		0 & \omega_2 & 0 & \cdot &\cdot & 0 \\
		\cdot &\cdot  & \cdot &\cdot  & \cdot &\cdot \\
		\cdot &\cdot  & \cdot &\cdot  & \cdot &\cdot \\
		0 & 0 & \cdot &\cdot & 0 & \omega_d
	\end{pmatrix},\]
	and let $f_1\in \cd.$ Then the following hold:
	
	\bei \item[(i)] If the vector $f_1\in \cd$ has
	a coordinate which is zero, then $\tnf$ is not
	a frame.
	\item[(ii)] If two or more of the roots
	$\omega_1, \omega_2, \dots, \omega_d$ are
	identical, then $\tnf$ is not a frame.
	\item[(iii)] If $\omega_1, \omega_2, \dots, \omega_d$ are distinct and all $d$ coordinates of
	the vector $f_1$ are nonzero, 
	then $\tnf$ is a cyclic frame for $\cd.$
	\item[(iv)] Under the assumptions in (iii), if at least one of the  $\omega_j, \, j=1, \dots,d,$ 
	is primitive, then $\tnf$ is a minimal cyclic frame.
	\eni 
	\el
	
	\bp First, writing 
	\bee \label{252002c} f_1= \begin{pmatrix}
		d_1 \\
		d_2 \\
		\vdots \\
		d_d
	\end{pmatrix}=  D \begin{pmatrix}
		1 \\
		1 \\
		\vdots \\
		1
	\end{pmatrix}, \, \, \, D:= \begin{pmatrix}
		d_1 & 0 & 0 & \cdot &\cdot & 0 \\
		0 & d_2 & 0 & \cdot &\cdot & 0 \\
		\cdot &\cdot  & \cdot &\cdot  & \cdot &\cdot \\
		\cdot &\cdot  & \cdot &\cdot  & \cdot &\cdot \\
		0 & 0 & \cdot &\cdot & 0 & d_d
	\end{pmatrix}, \ene we have
	\bee \label{252002a} T^{k-1} f_1 = D  T^{k-1} \begin{pmatrix}
		1 \\
		1 \\
		\vdots \\
		1
	\end{pmatrix}.\ene
	Thus it is clear that $\tnf$ can only be a frame if $D$ is invertible, which proves (i).
	
	Now, ignoring for the moment the matrix $D$ and
	looking  at the first $d$ vectors
	in \eqref{252002a},  we have
	
	\bee \label{252002b} \left\{ T^{k-1} \begin{pmatrix}
		1 \\
		1 \\
		\vdots \\
		1
	\end{pmatrix} \right\}_{k=1}^d= 
	\left\{ \begin{pmatrix} 1 \\ 1 \\ \cdot \\ \cdot \\ 1 \end{pmatrix}, \begin{pmatrix} \omega_1 \\ \omega_2 \\ \cdot \\ \cdot \\ \omega_d \end{pmatrix}, \begin{pmatrix} \omega_1^2 \\ \omega_2^2 \\ \cdot \\ \cdot \\ \omega_d^2 \end{pmatrix},
	\begin{pmatrix} \omega_1^3 \\ \omega_2^3 \\ \cdot \\ \cdot \\ \omega_d^{3} \end{pmatrix}, \cdots, \begin{pmatrix} \omega_1^{d-1} \\ \omega_2^{d-1} \\ \cdot \\ \cdot \\ \omega_d^{d-1} \end{pmatrix} \right\}.\ene
	
	The $d\times d$ matrix with the vectors  $  \{ T^{k-1}f_1 \}_{k=1}^d$ as columns is a Vandermonde matrix; if $\omega_j= \omega_\ell$ for
	some $j\neq \ell,$ the determinant vanishes,
	implying that the first $d$ vectors are not a basis for $\cd;$ hence, applying again \eqref{252002a} also the first $d$ vectors in $\tnf$ can not be a basis. Applying Lemma
	\ref{first_d_vector_basis} this proves (ii).
	
	On the other hand, in the case (iii) the
	Vandermonde matrix clearly has
	nonzero determinant; hence the vectors 
	in \eqref{252002b} are linearly independent, and therefore a basis for $\cd.$  Hence
	$\{ T^{k-1}f_1 \}_{k=1}^n$ is a  frame - and it is cyclic because $T^n=I,$  due to the choice of $\omega_1, \dots, \omega_d.$   If at least one of the $\omega_j, j=1, \dots,n,$ has order $n,$ clearly $T^\ell \neq I$ for all $\ell <n;$ this proves (iv).\ep

	Recall that a $d\times d$ matrix $T$ with $d$ distinct eigenvalues,
	$\omega_1, \omega_2, \dots, \omega_d,$ is diagonalizable, i.e., there exists an invertible $d\times d$ matrix $U$ such that
	\bee \label{251402a}  T  =U  \begin{pmatrix}
		\omega_1 & 0 & 0 & \cdot &\cdot & 0 \\
		0 & \omega_2 & 0 & \cdot &\cdot & 0 \\
		\cdot &\cdot  & \cdot &\cdot  & \cdot &\cdot \\
		\cdot &\cdot  & \cdot &\cdot  & \cdot &\cdot \\
		0 & 0 & \cdot &\cdot & 0 & \omega_d
	\end{pmatrix} U^{-1}. \ene Based on Lemma \ref{lem_construction_of_unitary_cyclic_dynamical_frames} we
	will now prove the following characterization of cyclic frames.

	\bt Fix $d,n\in \mn,$ and let $T$ denote an $d\times d$ matrix. Then the following hold:
	\bei \item[(i)] Assume that $T$  has $d$ distinct
	eigenvalues $\omega_1, \omega_2, \dots, \omega_d,$ each  being an $n$th root of unity. Let $f_1\in \cd$ denote an arbitrary vector with nonzero entries. Then, using the notation in \eqref{251402a} and letting
	$ \varphi:= Uf_1,$
	the sequence
	\bee \label{251402b} \{T^{n-1} \varphi\}_{k=1}^n\ene is
	a cyclic frame for $\cd.$ 
	\item[(ii)]  Assume that $\{T^{k-1}\varphi\}_{k=1}^n$ is
	a cyclic frame for some $\varphi \in \cd.$  Then $T$ is diagonalizable
	and has $d$ distinct eigenvalues
	$\omega_1, \omega_2, \dots, \omega_d,$ each  being an $n$th root of unity.
	Furthermore, with the notation in \eqref{251402a},
	the vector $\varphi$ has the form $\varphi=Uf_1$
	for a vector $f_1$ with nonzero coordinates.
	\eni

	\et
	
	\bp Under the assumptions in (i), direct calculation using \eqref{251402a}
	and applying again the notation in \eqref{252002c},
	
	\bes  T^{k-1} \varphi  & = &  \left[U  \begin{pmatrix}
		\omega_1 & 0 & 0 & \cdot &\cdot & 0 \\
		0 & \omega_2 & 0 & \cdot &\cdot & 0 \\
		\cdot &\cdot  & \cdot &\cdot  & \cdot &\cdot \\
		\cdot &\cdot  & \cdot &\cdot  & \cdot &\cdot \\
		0 & 0 & \cdot &\cdot & 0 & \omega_d
	\end{pmatrix} U^{-1}\right]^{k-1}  UD
	\begin{pmatrix}
		1 \\
		1 \\
		\vdots \\
		1 
	\end{pmatrix}    \\  & = & U    \begin{pmatrix}
		\omega_1 & 0 & 0 & \cdot &\cdot & 0 \\
		0 & \omega_2 & 0 & \cdot &\cdot & 0 \\
		\cdot &\cdot  & \cdot &\cdot  & \cdot &\cdot \\
		\cdot &\cdot  & \cdot &\cdot  & \cdot &\cdot \\
		0 & 0 & \cdot &\cdot & 0 & \omega_d
	\end{pmatrix}^{k-1}  D
	\begin{pmatrix}
		1 \\
		1 \\
		\vdots \\
		1 
	\end{pmatrix} \\ & = &  UD    \begin{pmatrix}
		\omega_1 & 0 & 0 & \cdot &\cdot & 0 \\
		0 & \omega_2 & 0 & \cdot &\cdot & 0 \\
		\cdot &\cdot  & \cdot &\cdot  & \cdot &\cdot \\
		\cdot &\cdot  & \cdot &\cdot  & \cdot &\cdot \\
		0 & 0 & \cdot &\cdot & 0 & \omega_d
	\end{pmatrix}^{k-1}  
	\begin{pmatrix}
		1 \\
		1 \\
		\vdots \\
		1 
	\end{pmatrix}.  \ens  Thus, 
	the sequence $ \{T^{n-1} \varphi\}_{k=1}^n$ is the image of the
	frame in Lemma \ref{lem_construction_of_unitary_cyclic_dynamical_frames} under the bijective map $UD,$ and hence a frame.  Applying
	\eqref{251402a} again, it follows immediately that the frame is
	cyclic.  
	
	In order to prove (ii), notice that if $T$ is not diagonalizable,
	its Jordan decomposition
	contains blocks of the form
	\bes   \begin{pmatrix}
		\lambda & 1 & 0 & \cdot &\cdot & 0 \\
		0 & \lambda & 1 & \cdot &\cdot & 0 \\
		\cdot &\cdot  & \cdot &\cdot  & \cdot &\cdot \\
		\cdot &\cdot  & \cdot &\cdot  & \cdot &\cdot \\ 	0 & 0 & \cdot &\cdot & \lambda & 1 \\
		0 & 0 & \cdot &\cdot & 0 & \lambda
	\end{pmatrix};
	\ens
	since such blocks are noncyclic, this excludes  that $T^n=I.$
	Thus $T$ is indeed diagonalizable. 
	The property $T^n=I$ immediately implies
	that all eigenvalues of $T$ are $n$th roots of unity; thus, we only have to prove that all
	eigenvalues have multiplicity one.  
	
	Considering again \eqref{251402a}, we
	have   \bes T^{k-1}  =U  \begin{pmatrix}
		\omega_1 & 0 & 0 & \cdot &\cdot & 0 \\
		0 & \omega_2 & 0 & \cdot &\cdot & 0 \\
		\cdot &\cdot  & \cdot &\cdot  & \cdot &\cdot \\
		\cdot &\cdot  & \cdot &\cdot  & \cdot &\cdot \\
		0 & 0 & \cdot &\cdot & 0 & \omega_d
	\end{pmatrix}^{k-1} U^{-1}; \ens assuming
	now that $\{T^{k-1}\varphi\}_{k=1}^n$ is 
	a frame, it follows that
	\bes \left\{    \begin{pmatrix}
		\omega_1 & 0 & 0 & \cdot &\cdot & 0 \\
		0 & \omega_2 & 0 & \cdot &\cdot & 0 \\
		\cdot &\cdot  & \cdot &\cdot  & \cdot &\cdot \\
		\cdot &\cdot  & \cdot &\cdot  & \cdot &\cdot \\
		0 & 0 & \cdot &\cdot & 0 & \omega_d
	\end{pmatrix}^{k-1} U^{-1}\varphi        \right\}
	\ens is a frame. By
	Lemma \ref {lem_construction_of_unitary_cyclic_dynamical_frames} (ii)+(i) this immediately implies
	that all the eigenvalues $\omega_j$ are distinct; furthermore, all coordinates of
	$U^{-1} \varphi$ are nonzero, meaning that
	$\varphi =   U U^{-1} \varphi =Uf_1,$ for a
	vector $f_1:= U^{-1} \varphi$ with nonzero coordinates. \ep 
	
	We now provide an alternative characterization
	of cyclic frames. 
	The  result is motivated by the infinite dimensional case  given in \cite[Theorem 2.1]{CHR2018}, but for technical reasons the proof (which we give in the appendix) is different. First, define
	the (cyclic) \textit{right shift operator} $R: \mc^n \to \mc^n$  by
	\begin{equation}\label{eq_Rz}
		Rx=R(x(1),	x(2),	\ldots,	x(n))
		=(x(n),	x(1),	\cdots,	x(n-1)),		
	\end{equation}	
	where $x=(x(1), \ldots, x(n)) \in \mc^n.$

	\bt \label{thm_cyclic_frame_invariant_kernel} Let $\fn$ be a frame in $\mc^d$ with synthesis operator $\Theta$. Then $\fn$ is a cyclic  frame if and only if $\ker(\Theta)$ is invariant under the  right shift operator $R.$
	\et

	The following result presents an alternative method for constructing cyclic frames; we  provide the proof in the appendix.
	
	\bt\label{thm_constructon_cyclic_dynamical_frame} 	Consider $d,n\in \mn$ with $n>d.$
	Assume that $a:=(a(1), \cdots, a(n))\in \mc^n$  has exactly $n-d$ non-zero coordinates. Let $c:= \check{a}$ denote
	the inverse 
	discrete Fourier transform of $a,$ and 
	consider the circulant matrix  
	\begin{equation}\label{eq_circulant}
		\mathcal{C}_n(c)= \begin{pmatrix}
			c(1) & c(n) & c(n-1) & \cdot & c(2) \\
			c(2) & c(1) & \cdot & \cdot & c(3) \\
			\cdot & \cdot & \cdot & \cdot & \cdot \\
			\cdot & \cdot & \cdot & \cdot & \cdot \\
			c(n) & c(n-1)  & \cdot &c(2) & c(1) 
		\end{pmatrix}_{n \times n}.
	\end{equation}	
	
	Define 
	$$M:=\{\mathcal{C}_n(c)z \ : \ z \in \C^n \}.$$
	Then, choosing the set $\{v_1, v_2, \cdots, v_d\}$ as a basis for the orthogonal complement $M^\perp$ of $M$ in $\mc^n$ and letting $V:=[v_1, v_2, \cdots, v_d]$,  the columns of the transposed matrix $V^t$ forms a cyclic  frame for $\cd$.
	\et
	
	The following application of Theorem~\ref{thm_constructon_cyclic_dynamical_frame} uses an auxiliary result from the appendix, see Lemma~\ref{lem_eigenvalues_circulant}.
	\bex	\label{ex_circulant}
		Suppose $d=3$, $n=4$ and $a=\begin{pmatrix}
			0 &1 &0 &0 
		\end{pmatrix} \in \C^4$.	Define $\mathcal{C}_4(c)$ to be  a circulant  matrix  by taking $c=\check{a}$ as follows 
		$$\mathcal{C}_4(c)=\begin{pmatrix}
			1/4 &-i/4 &-1/4 &i/4\\
			i/4&1/4 &-i/4 &-1/4\\
			-1/4	&i/4 &1/4 &-i/4\\
			-i/4	&-1/4 &i/4 &1/4
		\end{pmatrix}.$$
		Define $M:= \mbox{range}(\mathcal{C}_4(c))$.  Using Lemma~\ref{lem_eigenvalues_circulant},  only one eigenvalue of $\mathcal{C}$ is non-zero and hence dim$(M)=1$. Therefore,  we have   $V^t= \begin{pmatrix}
			-i &1 & 0 &0\\
			1 &0 &1 &0\\
			i &0 &0 &1
		\end{pmatrix}$. Thus columns of $V^t$ forms a cyclic  frame $\mc^3$, using Theorem \ref{thm_constructon_cyclic_dynamical_frame}.  Moreover, if we consider
		$$\{T^{k-1} f_1\}_{k=1}^4=\left\{\begin{pmatrix}
			-i \\1 \\ i
		\end{pmatrix}, \begin{pmatrix}
			1 \\0 \\0
		\end{pmatrix}, \begin{pmatrix}
			0\\	1 \\0
		\end{pmatrix},\begin{pmatrix}
			0\\0 \\1
		\end{pmatrix} \right\},$$
		then
		$$T=\begin{pmatrix}
			0 &0 &-i\\
			1 &0 &1\\
			0 &1 &i
		\end{pmatrix} \mbox{ and }f_1=\begin{pmatrix}
			-i \\1\\i
		\end{pmatrix}.$$
		\enx

	We end this section with a number of observations that are needed in Section \ref{251302c}. We first prove a relationship between the frame bounds for a cyclic frame $\tnf$ and the norm of the operators $T$ and $T^{-1}.$ 
	
	\bl \label{250603b} Let $\tnf$ be a cyclic frame for $\cd,$ with frame bounds $A,B.$ Then
	\bes 1 \le ||T|| \le \sqrt{\frac{B}{A}} \, \, \, \mbox{and}  \, \, \,  1 \le ||T^{-1}|| \le \sqrt{\frac{B}{A}}. \ens \el 
	\bp The proof of the estimate on $||T||$ follows
	the same steps as the proof for frames of the form $\{T^k \varphi\}_{k=-\infty}^\infty$ in an
	infinite-dimensional Hilbert space, see
	\cite[Theorem~2.3]{olemmaarzieh-2}. Since 
	$\{ \left(T^{-1}\right)^{k-1}[ T^{n-1}f_1]\}_{k=1}^n$ consists of precisely
	the same vectors as $\tnf$ (just ``in the opposite order), the result about the norm of
	$T^{-1}$ is an immediate consequence. \ep
	
	It is well-known that the canonical dual frame of a dynamical frame $\tnf$ is again a dynamical frame; indeed, denoting the frame operator for $\tnf$ by $S,$ we have that
	\bee \label{251203a} \{\si T^{k-1}f_1\}_{k=1}^n=\{(S^{-1} T S)^{k-1}S^{-1}f_1 \}_{k=1}^{n}. \ene  
	For a cyclic frame $\tnf,$  the operator
	$\si TS$ equals the inverse of the adjoint
	$T^*:$ 
	
	\bl \label{prop_canonical_dynamical_dual}
	Let $\tnf$ be a cyclic frame with frame operator $S$. Then  $S^{-1} T S=(T^\ast)^{-1}$.
	\el
	
	\bp Using that $T^n=I,$ a direct calculation shows that for any $f \in \cd,$
	\bee \label{eq_TSTstar=T} TSf= T\sum_{k=1}^{n}\langle  f, T^{k-1}f_1\rangle T^{k-1}f_1= \sum_{k=1}^{n}\langle (T^\ast)^{-1} f, T^{k}f_1\rangle T^{k}f_1=S(T^\ast)^{-1} f. \ene
	Hence $ T S=S(T^\ast)^{-1},$ and the result follows.
	\ep

	\section{Tight cyclic frames and erasures} \label{251302c}

	In this section we  show that tight cyclic frames have desirable properties in the context of erasure problems.
	First, recall that a frame $\fn$ is called an {\it equal-norm} or \textit{uniform frame} if $\|f_k\|=c$ for each $k \in \{1,\ldots,n\}$. A frame $\fn$ for $\cd$ is called {\it equiangular} or \textit{$2$-uniform frame} if it is uniform  and $|\la f_{k_1}, f_{k_2} \ra| =c$ for all $k_1 \neq k_2$. These frames play an important role in the case of erasures. For example, in \cite{HRBPVI} it is proved that a Parseval frame $\ftn$ is optimal for $1$-erasures if and only if it is a uniform frame; 
	further,  if a frame $\ftn$ is optimal  for $1$-erasures and  it is equiangular, then it is optimal for $ 2$ erasures. We refer to \cite{HRBPVI}  for the exact definition of these concepts.
	
	In order to prove our main results about tight cyclic frames, we need the next result, which is a direct consequence of Lemma \ref{prop_canonical_dynamical_dual}.

	\bl \label{prop_tight_cyclic_unitary_operator}
	If $\tnf$ is tight cyclic frame for $\cd$, then $T$ is  unitary. 
	\el

	We now prove that tight cyclic frames are
	uniform, and provide a characterization of equiangularity.

	\bpr \label{250603a} Let $\tnf$ be a tight cyclic frame  $\cd.$ Then the following hold: 
	\bei\item[(i)] $\tnf$ is uniform.
	\item[(ii)]   $\tnf$ is equiangular if and only if $| \la T^\ell f_1, f_1 \ra|$ is
	constant for $\ell \in \{1,2,\ldots, n-1\}.$\eni 
	\epr
	
	\bp Concerning (i), it follows directly from Lemma \ref{prop_tight_cyclic_unitary_operator} that the operator $T$ is unitary; thus 
	$||T^{k-1}f_1|| = ||f_1||.$  Concerning (ii), 
	since the operator $T$ is unitary 
	\bes 
	|\la T^{i}f_1, T^{j}f_1 \ra |=| \la T^{i-j} f_1 , f_1 \ra|, \, i.j=1, \dots,n; \ens thus
	the result follows from (i). 
	\ep

	Let us phrase a concrete version of the results in Lemma \ref{prop_tight_cyclic_unitary_operator} and    Proposition \ref{250603a} for
	the canonical tight frame associated with any cyclic frame:

	\bc \label{cor_equiangular_characterization_canonical_tight_frame} Given any
	cyclic frame $\tnf$ with frame operator $S,$ the following holds: 
	\bei
	\item[(i)]  The sequence $ \left\{ \left( S^{-1/2}TS^{1/2}\right)^{k-1} S^{-1/2}f_1\right\}_{k=1}^n $ is a cyclic, tight, and equal-norm frame. 
	
	\item[(ii)]  The operator $S^{-1/2}TS^{1/2} $ is unitary.
	
	\item[(iii)]  The sequence  
	$ \left\{ \left( S^{-1/2}TS^{1/2}\right)^{k-1} S^{-1/2}f_1\right\}_{k=1}^n$ is equiangular if and only if $  | \la T^\ell f_1, S^{-1}f_1 \ra|    $ is constant for
	$\ell \in \{1,2,\ldots, n-1\}.$ 
	\eni
	\ec
	
\bp
It is well known  that if  $\tnf$ is frame, then $\left\{S^{-1/2} (T^{k-1} f_1)\right\}_{k=1}^n$ is a tight frame.  Moreover, observe that 
\bes
\left\{ \left( S^{-1/2}TS^{1/2}\right)^{k-1} S^{-1/2}f_1\right\}_{k=1}^n =\left\{S^{-1/2} (T^{k-1} f_1)\right\}_{k=1}^n.
\ens
Therefore, $ \left\{ \left( S^{-1/2}TS^{1/2}\right)^{k-1} S^{-1/2}f_1\right\}_{k=1}^n $ is a tight frame for $\cd$. Next, we prove that this frame is cyclic. A direct calculation gives 
\bes
\left( S^{-1/2}TS^{1/2}\right)^{n}=S^{-1/2}T^{n}S^{1/2}= S^{-1/2}S^{1/2}=I,
\ens
where we have used that $T^n=I$. Hence the sequence is indeed a cyclic tight frame for $\cd$. Now the result follows from Lemma \ref{prop_tight_cyclic_unitary_operator} and Proposition \ref{250603a}.
\ep
	
	Lemma \ref{251302g} and Corollary \ref{cor_equiangular_characterization_canonical_tight_frame} lead to the following  construction of an equiangular cyclic frame $\tnf$ for $\cd$ with $n=d+1$. 
	\bpr \label{252003a} Let  $\fd$ be the standard orthnormal basis for $\cd,$ and 
	consider the cyclic frame \(\{T^{k-1}f_1\}_{k=1}^{d+1}=\{f_k\}_{k=1}^{d+1}\), where $f_{d+1}=-f_1-f_2- \cdots - f_d.$  Then, letting $S$ denote the frame operator, the sequence  
	$$ \left\{ \left( S^{-1/2}TS^{1/2}\right)^{k-1} S^{-1/2}f_1\right\}_{k=1}^{d+1}$$
	is a cyclic, tight, equal-norm and equiangular frame for $\cd.$  
	\epr
	
	\bp
	A direct calculation shows that
	the frame operator $S$ 
	and its inverse $S^{-1}$ are given by 
	\bee \label{eq_S_and_its_inverse}
	S=	\begin{pmatrix}
		2& 1 & 1 & \cdot &\cdot & 1 \\
		1 & 2 & 1 & \cdot &\cdot & 1 \\
		1 & 1 & 2 & \cdot &\cdot & 1 \\
		\cdot &\cdot  & \cdot &\cdot  & \cdot &\cdot \\
		\cdot &\cdot  & \cdot &\cdot  & \cdot &\cdot \\
		1 & 1 & \cdot &\cdot & 1 & 2
	\end{pmatrix}_{d \times d} \mbox{ and } S^{-1}=\begin{pmatrix}
		\frac{d}{d+1} & -\frac{1}{d+1} & -\frac{1}{d+1} & \cdot &\cdot & -\frac{1}{d+1} \\
		-\frac{1}{d+1} & \frac{d}{d+1} & -\frac{1}{d+1} & \cdot &\cdot & -\frac{1}{d+1} \\
		-\frac{1}{d+1} &	-\frac{1}{d+1} & \frac{d}{d+1} &  \cdot &\cdot & -\frac{1}{d+1} \\
		\cdot &\cdot  & \cdot &\cdot  & \cdot &\cdot \\
		\cdot &\cdot  & \cdot &\cdot  & \cdot &\cdot \\
		-\frac{1}{d+1} & -\frac{1}{d+1} & \cdot &\cdot & -\frac{1}{d+1} & \frac{d}{d+1}
	\end{pmatrix}_{d \times d}.
	\ene
	Therefore 
	\bes S^{-1}f_1=\begin{bmatrix}
		\frac{d}{d+1}\\
		-\frac{1}{d+1}\\
		\cdot\\
		\cdot\\
		-\frac{1}{d+1}
	\end{bmatrix}=f_1+ \begin{bmatrix}
		-\frac{1}{d+1}\\
		-\frac{1}{d+1}\\
		\cdot\\
		\cdot\\
		-\frac{1}{d+1}
	\end{bmatrix}=f_1 -\frac{1}{d+1}(f_1+f_2+\ldots+f_d).\ens
	
	Thus, for $\ell \in \{1,\ldots,d-1\}$, we have 
	\bes |\la T^\ell f_1, S^{-1}f_1 \ra|=  \left|\frac{-1}{d+1}\la f_{\ell+1},  f_{\ell+1} \ra \right|= \frac{1}{d+1},\ens
	and with $\ell=d,$ 
	\bes
	|\la T^d f_1, S^{-1}f_1 \ra|&=&\left|-\la f_1, f_1 \ra+ \frac{1}{d+1}(\la f_1, f_1\ra+ \ldots+\la f_d, f_d\ra)  \right|\\
	&=& \left| -1+\frac{d}{d+1} \right|= \frac{1}{d+1}. 
	\ens	
	Thus,  $|\la T^\ell f_1, S^{-1}f_1 \ra|= \frac{1}{d+1}$ for each $\ell \in \{1,\ldots,d\}$.  Now the result follows from  Corollary \ref{cor_equiangular_characterization_canonical_tight_frame} (i)+(iii). 
	\ep
	
	The following example demonstrates that applying the construction in Proposition 4.4 in the case of $\mr^2$ results in a frame which is a rotation of the Mercedes-Benz frame.
	\bex \label{252103c}
	Consider the cyclic frame $\{T^{k-1}f_1\}_{k=1}^3=\{f_k\}_{k=1}^3$ in $\mr^2$ with frame operator $S$, where $\{f_1, f_2\}$ are standard orthonormal basis for $\mr^2$, and $f_3=-f_1-f_2$.  The  operators $S$ and $S^{-1}$ are the same as those defined in equation  (\ref{eq_S_and_its_inverse}) for $d=2$. Thus, a direct calculation gives the operator $S^{1/2}$ 
	and its inverse $S^{-1/2}$ as follows: 
	\bes
	S^{1/2}=	\begin{pmatrix}
		\frac{\sqrt{3}-1}{2} & \frac{\sqrt{3}+1}{2}\\
		\frac{\sqrt{3}+1}{2} & \frac{\sqrt{3}-1}{2}
	\end{pmatrix} \mbox{ and } S^{-1/2}=\begin{pmatrix}
		\frac{1-\sqrt{3}}{2 \sqrt{3}} & \frac{1+\sqrt{3}}{2 \sqrt{3}}\\
		\frac{1+\sqrt{3}}{2 \sqrt{3}} & \frac{1-\sqrt{3}}{2 \sqrt{3}}
	\end{pmatrix}.
	\ens
	
	Hence,  by Proposition 4.4, the sequence  \bes\left\{ \left( S^{-1/2}TS^{1/2}\right)^{k-1} S^{-1/2}f_1\right\}_{k=1}^{3}= \left\{ \begin{pmatrix}
		\frac{1-\sqrt{3}}{2 \sqrt{3}} \\
		\frac{1+\sqrt{3}}{2 \sqrt{3}}
	\end{pmatrix},  \begin{pmatrix}
		\frac{1+\sqrt{3}}{2 \sqrt{3}} \\
		\frac{1-\sqrt{3}}{2 \sqrt{3}}
	\end{pmatrix}, \begin{pmatrix}
		-\frac{1}{ \sqrt{3}} \\
		-\frac{1}{ \sqrt{3}}
	\end{pmatrix} \right\}\ens
	is a tight  equiangular frame for $\mr^2$. An elementary calculation
	reveals that the frame is a rotation 
	with $\frac{-\pi}{12}$ rad of the Mercedes Benz frame in Figure  1.\enx

	\bex \label{ex_equiangular_cyclic_frame_for_C^3}
	Recall the cyclic  frame constructed in Example \ref{ex_circulant}
	$$\{T^{k-1} f_1\}_{k=1}^4=\left\{\begin{pmatrix}
		-i \\1 \\ i
	\end{pmatrix}, \begin{pmatrix}
		1 \\0 \\0
	\end{pmatrix}, \begin{pmatrix}
		0\\	1 \\0
	\end{pmatrix},\begin{pmatrix}
		0\\0 \\1
	\end{pmatrix} \right\},$$
	for $\mc^3$, where $$T=\begin{pmatrix}
		0 &0 &-i\\
		1 &0 &1\\
		0 &1 &i
	\end{pmatrix} \mbox{ and }f_1=\begin{pmatrix}
		-i \\1\\i
	\end{pmatrix}.$$  Let $S$ be the frame operator of $\{T^{k-1} f_1\}_{k=1}^4$. Then
	$$S= \begin{pmatrix}
		2 &-i &-1\\
		i &2 &-i\\
		-1 &i &2
	\end{pmatrix} \mbox{ and } S^{-1}=\begin{pmatrix}
		\frac{3}{4} &\frac{i}{4} &\frac{1}{4}\\
		-\frac{i}{4} &\frac{3}{4} &\frac{i}{4}\\
		\frac{1}{4} &-\frac{i}{4} &\frac{3}{4}
	\end{pmatrix}.$$
	Now $S^{-1}f_1=\begin{pmatrix}
		-\frac{i}{4} \\ \frac{1}{4} \\ \frac{i}{4}
	\end{pmatrix}$. A simple calculations shows that for each $\ell \in \{1,2,3\}$ we have 
	$$|\la T^\ell f_1, S^{-1}f_1 \ra|=\frac{1}{4}.$$
	Hence, by Corollary \ref{cor_equiangular_characterization_canonical_tight_frame} (iii), the frame
	$\left\{ \left( S^{-1/2}TS^{1/2}\right)^{k-1} S^{-1/2}f_1\right\}_{k=1}^4$ is equiangular.  In this case 
	$$\ S^{\frac{1}{2}}=\begin{pmatrix}
		\frac{4}{3} &-\frac{i}{3} &-\frac{1}{3}\\
		\frac{i}{3} &\frac{4}{3} &-\frac{i}{3}\\
		-\frac{1}{3} &\frac{i}{3} &\frac{4}{3}
	\end{pmatrix} \mbox{ and } S^{-\frac{1}{2}}=\begin{pmatrix}
		\frac{5}{6} &\frac{i}{6} &\frac{1}{6}\\
		-\frac{i}{6} &\frac{5}{6} &\frac{i}{6}\\
		\frac{1}{6} &-\frac{i}{6} &\frac{5}{6}
	\end{pmatrix}.$$
	Therefore,
	$$\{(S^{-\frac{1}{2}}TS^{\frac{1}{2}})^{k-1} S^{-\frac{1}{2}}f_1\}_{k=1}^{4}=\left\{\begin{pmatrix}
		-\frac{i}{2} \\
		\frac{1}{2}\\
		\frac{i}{2}
	\end{pmatrix}, \begin{pmatrix}
		\frac{5}{6}\\
		-\frac{i}{6}\\
		\frac{1}{6}
	\end{pmatrix}, \begin{pmatrix}
		\frac{i}{6}\\
		\frac{5}{6}\\
		-\frac{i}{6}
	\end{pmatrix}, \begin{pmatrix}
		\frac{1}{6}\\
		\frac{i}{6}\\
		\frac{5}{6}
	\end{pmatrix} \right\}.$$
	\enx

	\appendix \section{Appendix: Proof of Thm.  \ref{thm_cyclic_frame_invariant_kernel}  and  Thm. \ref{thm_constructon_cyclic_dynamical_frame} }

	\noindent{\bf Proof of Theorem \ref{thm_cyclic_frame_invariant_kernel}: }
	If $\ftkn$ is a cyclic frame for $\cd,$
	a direct calculations shows that if  $x \in \ker(\Theta), $ then also $R(x) \in \ker(\Theta) $.  Conversely assume that  $R[\ker (\Theta)] \subseteq \ker(\Theta)$. Define the operator $T: \cd \to \cd$ by
	$$T\left(\sum\limits_{k=1}^{n} \alpha(k) f_k\right)=\alpha(n)f_1+ \sum_{k=1}^{n-1} \alpha
	(k)f_{k+1}.$$
	First, we  show that $T$ is well-defined. Suppose $\sum\limits_{k=1}^{n} \alpha(k) f_k=\sum\limits_{k=1}^{n} \beta(k) f_k; then, $ letting   $\gamma(k)=\alpha(k)-\beta(k)$ for $1 \leq k \leq n$, we have that $\gamma \in \ker(\Theta)$ and thus  $R(\gamma) \in \ker(\Theta);$ that is, 
	$$0= \gamma(n)f_1+\sum\limits_{k=2}^{n}  \gamma(k-1) f_{k}=(\alpha(n)-\beta(n))f_{1}+\sum_{k=2}^{n-1} (\alpha(k)-\beta(k))f_{k+1}.$$
	It immediately follows that  
	$\alpha(n)f_1+\sum_{k=1}^{n-1} \alpha(k) f_{k+1}=\beta(n)f_1+\sum_{k=1}^{n-1} \beta(k) f_{k+1},$ and hence that
	$T$ is well defined.  It is now easy to verify that $T^nf_1=f_1,$ implying that  $\fn$ is a cyclic frame.
	\ep
	
	In order to prove Theorem \ref{thm_constructon_cyclic_dynamical_frame}
	we first state a number of auxiliary results
	about circulant matrices. For these results
	we consider an arbitrary circulant matrix 
	$\mathcal{C}_n(c),$  i.e., we do not
	assume that the vector $c$ has the
	special form in Theorem \ref{thm_constructon_cyclic_dynamical_frame}.

	\bl \label{lem_circulant_Invariant}
	If $\mathcal{C}_n(c)$  is a circulant matrix, then $\mbox{range}(\mathcal{C}_{n}(c))$  is invariant under the right shift operator $R$. 
	\el  
	\bp
	Let $y \in \mbox{range}(\mathcal{C}_{n}(c));$ then, by the definition of a circulant matrix, as defined in (\ref{eq_circulant}), 	$ y=x(1)c+x(2)Rc+\ldots+x(n)R^{n-1}c $
	for some $x=(x(1), \ldots, x(n)) \in \mcn.$
	Therefore
	\bes Ry=x(1)Rc+x(2)R^2c+\ldots+x(n-1)R^{n-1}c+x(n)c \in  \mbox{range}(\mathcal{C}_{n}(c)), \ens
	as claimed.
	\ep
	The following result and its proof can be found in \cite[Corollary 3.33]{PPSTBook}.
	
	\bl \label{lem_eigenvalues_circulant}
	If $\mathcal{C}_n(c)$  is a circulant matrix, then the set of eigenvalues  is given by $ \{\widehat{c}(1), \cdots, \widehat{c}(n)\}$, where $\widehat{c}$ denotes the discrete Fourier transform of $c$ .
	\el
	
	\vn{\bf Proof of Theorem~\ref{thm_constructon_cyclic_dynamical_frame}:}
	By Lemma \ref{lem_circulant_Invariant}, $M$ is invariant under right shift operator $R.$ Using Lemma $\ref{lem_eigenvalues_circulant}$, the set of eigenvalues  $\mathcal{C}_n(c)$ is $ \{\widehat{c}(1), \cdots, \widehat{c}(n)\}$.  Since $\widehat{c}(k)\neq 0$ for $n-d$ values of $k \in \{1,\cdots,n\}$, it follows that $\dim(M)=n-d$. Now,  the matrix $V=[v_1, v_2, \cdots, v_d],$ has rank $d$. 
	Then the transposed matrix $V^t$  is of order  $d \times n$ having the property $$V^tx=0 \mbox{ for all } x \in M.$$ This implies that $\ker(V^t) \supset M$. Indeed  $\ker(V^t)=M$, follows by observing that $\dim(M)=n-d$, and  the dimension of  $\ker(V^t) $ is $n-d$ since the rank of  $V^t$ is  $d$. Thus, $V^t$ is a matrix of size $d \times n$ with kernel $M$, which is invariant under the right shift operator. Hence by Theorem \ref{thm_cyclic_frame_invariant_kernel}, the columns of $V^t$ forms a cyclic frame for $\cd$.
	\ep
	
	For the sake of convenient application, we 
	note that Theorem \ref{thm_constructon_cyclic_dynamical_frame}
	can be phrased as an algorithm, as follows.

	\vspace{.5cm}
	\hrule
	\vspace{.1cm}
	\textbf{Algorithm 1:} Construction of cyclic  frame
	\vspace{.1cm}
	\hrule
	\vspace{.1cm}
	\noindent\textbf{Input:} \\
	$d,n \in \mn$ \\
	A vector $a \in \mc^n$ has $n-d$ entry non-zero.\\
	\textbf{Output:} $V^t$ is a matrix of order  $d \times n$ having kernel $M$. \\
	Set $c=\check{a}$\\
	Set $\mathcal{C}_n(c)$  as a circulant matrix having first column $c$.\\
	Set $M=\mbox{Range}(\mathcal{C}_n(c))$\\
	Set $M^\perp$ orthogonal complement of $M$\\
	Set $V$ is a matrix of order $n \times d$ such that columns vector form a basis of $M^\perp$\\
	Set $V^t$ transpose of $V$\\
	The final output is $V^t$. Then by Theorem \ref{thm_constructon_cyclic_dynamical_frame}, the columns vector of  $V^t$ forms a cyclic  frame for $\C^d$. 
	\vspace{.1cm}
	\hrule
	\vspace{.3cm}

	\section*{Acknowledgements}
 The authors would like to express their sincere thanks to anonymous reviewers for carefully reading the manuscript and providing us with  valuable comments and suggestions, which improved the presentation of this manuscript. The second and third authors are also grateful to Lokendra Kumar for fruitful discussions at the beginning of this work.
	
		\vspace{.3cm}
		
		\section*{Funding} The research of Navneet Redhu was supported by a research grant from CSIR, New Delhi [09/1022(0099)/2020-EMR-I].  Niraj K. Shukla's research was supported by the DST-SERB Project [MTR/2022/000176].

	{\bf \vspace{.1in} \noindent Ole Christensen\\
		Department of Applied Mathematics and Computer Science\\
		Technical University of Denmark,
		Building 303,
		2800 Lyngby, Denmark\\
		Email: ochr@dtu.dk}
	
	{\bf \vspace{.1in} \noindent Navneet Redhu\\
		Department of Mathematics\\
		Indian Institute of Technology Indore\\
		Simrol, Khandwa Road,
		Indore-453 552, Madhya Pradesh, India \\
		Email: phd2001141011@iiti.ac.in}
	
	{\bf \vspace{.1in} \noindent Niraj K. Shukla\\
		Department of Mathematics\\
		Indian Institute of Technology Indore\\
		Simrol, Khandwa Road,
		Indore-453 552, Madhya Pradesh, India \\
		Email: nirajshukla@iiti.ac.in}

\end{document}